\documentclass[leqno,12pt]{article}

\usepackage{amssymb}
\usepackage{euscript}
\usepackage[dvips]{graphicx}
\usepackage{flafter}
\usepackage{graphicx}
\usepackage{pstricks}

\evensidemargin\oddsidemargin

\begin{document}

\baselineskip=18pt
\setcounter{page}{1}

\renewcommand{\theequation}{\thesection.\arabic{equation}}
\newtheorem{theorem}{Theorem}[section]
\newtheorem{lemma}[theorem]{Lemma}
\newtheorem{proposition}[theorem]{Proposition}
\newtheorem{corollary}[theorem]{Corollary}
\newtheorem{remark}[theorem]{Remark}
\newtheorem{fact}[theorem]{Fact}
\newtheorem{problem}[theorem]{Problem}

\newcommand{\eqnsection}{
\renewcommand{\theequation}{\thesection.\arabic{equation}}
    \makeatletter
    \csname  @addtoreset\endcsname{equation}{section}
    \makeatother}
\eqnsection

\def\r{{\mathbb R}}
\def\e{{\mathbb E}}
\def\p{{\mathbb P}}
\def\P{{\bf P}}
\def\E{{\bf E}}
\def\Q{{\bf Q}}
\def\z{{\mathbb Z}}
\def\T{{\mathbb T}}
\def\G{{\mathbb G}}

\def\ee{\mathrm{e}}
\def\d{\, \mathrm{d}}
\def\deg{{\mathrm b}}




\vglue50pt

\centerline{\Large\bf Slow movement of random walk}

\medskip

\centerline{\Large\bf in random environment on a regular tree}

\bigskip
\bigskip

\centerline{by}

\medskip

\centerline{Yueyun Hu $\;$and$\;$ Zhan Shi}

\medskip

\centerline{\it Universit\'e Paris XIII \& Universit\'e Paris VI}

\bigskip
\bigskip

{\leftskip=2truecm
\rightskip=2truecm
\baselineskip=15pt
\small

\noindent{\slshape\bfseries Summary.} We consider a recurrent random walk in random environment on a regular tree. Under suitable general assumptions upon the distribution of the environment, we show that the walk exhibits an unusual slow movement: the order of magnitude of the walk in the first $n$ steps is $(\log n)^3$.

\bigskip


\bigskip


} 

\bigskip
\bigskip

\section{Introduction}
   \label{s:intro}

Let $\T$ be a rooted $\deg$-ary tree, with $\deg\ge 2$. Let $\omega:= (\omega(x,y), \, x,y\in \T)$ be a collection of non-negative random variables such that $\sum_{y\in \T} \omega(x,y)=1$ for any $x\in \T$. Given $\omega$, we define a Markov chain $X:= (X_n, \, n\ge 0)$ on $\T$ with $X_0 =e$ and
$$
P_\omega(X_{n+1}= y \, | \, X_n =x) = \omega(x, y) .
$$

\noindent The process $X$ is called random walk in random environment (or simply RWRE) on $\T$. (By informally taking $\deg=1$, $X$ would become a usual RWRE on the half-line $\z_+$.)

We refer to page 106 of Pemantle and Peres~\cite{pemantle-peres1} for a list of motivations to study tree-valued RWRE. For a close relation between tree-valued RWRE and Mandelbrot's multiplicative cascades, see Menshikov and Petritis~\cite{menshikov-petritis}.

We use $\P$ to denote the law of $\omega$, and the semi-product measure $\p (\cdot) := \int P_\omega (\cdot) \P(\! \d \omega)$ to denote the distribution upon average over the environment.

Some basic notation of the tree is in order. Let $e$ denote the root of $\T$. For any vertex $x\in \T \backslash \{ e\}$, let ${\buildrel \leftarrow \over x}$ denote the parent of $x$. As such, each vertex $x\in \T \backslash \{ e\}$ has one parent ${\buildrel \leftarrow \over x}$ and $\deg$ children, whereas the root $e$ has $\deg$ children but no parent. For any $x\in \T$, we use $|x|$ to denote the distance between $x$ and the root $e$: thus $|e|=0$, and $|x| = |{\buildrel \leftarrow \over x}| +1$.

We define
\begin{equation}
    A(x) := {\omega({\buildrel \leftarrow \over x},
    x) \over \omega({\buildrel \leftarrow \over x},
    {\buildrel \Leftarrow \over x})} , \qquad
    x\in \T, \; |x|\ge 2 ,
    \label{A}
\end{equation}

\noindent where ${\buildrel \Leftarrow \over x}$ denotes the parent of ${\buildrel \leftarrow \over x}$.

Following Lyons and Pemantle~\cite{lyons-pemantle}, we assume throughout the paper that $(\omega(x,\bullet))_{x\in \T\backslash \{ e\} }$ is a family of i.i.d.\ {\it non-degenerate} random vectors and that $(A(x), \; x\in \T, \; |x|\ge 2)$ are identically distributed. We also assume the existence of $\varepsilon_0>0$ such that $\omega(x,y) \ge \varepsilon_0$ if either $x= {\buildrel \leftarrow \over y}$ or $y= {\buildrel \leftarrow \over x}$, and $\omega(x,y) =0$ otherwise; in words, $(X_n)$ is a nearest-neighbour walk, satisfying an ellipticity condition.

Let $A$ denote a generic random variable having the common distribution of $A(x)$ (for $|x| \ge 2$) defined in (\ref{A}). Let
\begin{equation}
    p := \inf_{t\in [0,1]} \E (A^t) .
    \label{p}
\end{equation}

\noindent An important criterion of Lyons and Pemantle~\cite{lyons-pemantle} says that with $\p$-probability one, the walk $(X_n)$ is recurrent or transient, according to whether $p\le {1\over \deg}$ or $p>{1\over \deg}$. It is, moreover, positive recurrent if $p<{1\over \deg}$. Later, Menshikov and Petritis~\cite{menshikov-petritis} proved that the walk is null recurrent if $p={1\over \deg}$.

Throughout the paper, we write
$$
X_n^* := \max_{0\le k\le n} |X_k|, \qquad n\ge 0.
$$

In the positive recurrent case $p<{1\over \deg}$, ${X_n^*\over \log n}$ converges $\p$-almost surely to a constant $c\in (0,\, \infty)$ whose value is known, see \cite{yztree}.

The null recurrent case $p={1\over \deg}$ is more interesting. It turns out that the behaviour of the walk depends also on the sign of $\psi'(1)$, where
\begin{equation}
    \psi(t) := \log \E \left( A^t \right), \qquad t
    \ge 0.
    \label{psi}
\end{equation}

\noindent In \cite{yztree}, we proved that if $p={1\over \deg}$ and $\psi'(1)<0$, then
\begin{equation}
    \lim_{n\to \infty} \, {\log X_n^* \over \log n}
    = 1- {1\over \min\{ \kappa, 2\} } \; ,
    \qquad \hbox{\rm $\p$-a.s.},
    \label{yztree}
\end{equation}

\noindent where $\kappa := \inf \{ t>1: \; \E(A^t) = {1\over \deg} \} \in (1, \, \infty]$, with $\inf \emptyset :=\infty$.

The delicate case $p={1\over \deg}$ and $\psi'(1)\ge 0$ was left open, and is studied in the present paper. See Figure 1.

\bigskip

\begin{figure}[h]
\centerline{%
\hbox{%
\psset{unit=4.5mm} \pspicture(0,0)(39,10)
 \Cartesian(4.5mm,4.5mm)
\psline{->}(8,0)(8,10) \rput[rb](7.8,6.5){\small $0$}
\psline{->}(8,7)(16, 7) \rput[rb](16.5,6.5){\small $t$}
\rput[rb](9,10.2){{\small $\psi(t)$}}
\psbezier[linewidth=1.1pt]{-}(8,7)(9,4)(12,-0.6)(14, 4)
\psline[linewidth=0.5pt,linestyle=dashed]{-}(8,2)(16,2)
\rput[rb](7.7,1.6){\small $ \log {1\over \deg} $}
\psline[linewidth=0.5pt,linestyle=dashed]{-}(12,7)(12,2)
\rput[rb](12.2,7.2){\small $1$}
\psline{->}(24,0)(24,10) \rput[rb](23.8,6.5){\small $0$}
\psline{->}(24,7)(32, 7) \rput[rb](32.5,6.5){\small $t$}
\rput[rb](25,10.2){{\small $\psi(t)$}}
\psbezier[linewidth=1.1pt]{-}(24,7)(26.3,-2.2)(28,4)(28.6,  5)
\psline[linewidth=0.5pt,linestyle=dashed]{-}(24,2)(32,2)
\rput[rb](23.7,1.6){\small $ \log {1\over \deg} $}
\psline[linewidth=0.5pt,linestyle=dashed]{-}(26.5,7)(26.5, 2)
\rput[rb](26.8,7.2){\small $ \theta $}
\psline[linewidth=0.5pt,linestyle=dashed]{-}(28,7)(28,3.8)
\rput[rb](28.2,7.2){\small $1$}
\endpspicture}}
\caption{ \textsf{Case  $\psi'(1)= 0$ and case $\psi'(1)>0$ with
 $\theta $   defined in (\ref{gamma}).}}
\end{figure}

We will see in Remark \ref{r:V*} that the case $\psi'(1)>0$ boils down to the case $\psi'(1)=0$ via a simple transformation of the distribution of the random environment. As is pointed out by Biggins and Kyprianou~\cite{biggins-kyprianou} in the study of Mandelbrot's multiplicative cascades, the case $\psi'(1)=0$ is likely to be ``both subtle and important".

The following theorem reveals an unusual slow regime for the walk.

\medskip

\begin{theorem}
 \label{t:main}
 If $p={1\over \deg}$ and if $\psi'(1)\ge 0$, then
 there exist constants $0<c_1 \le c_2 <\infty$ such
 that
 \begin{equation}
     c_1 \le \liminf_{n\to \infty} {X_n^*\over
     (\log n)^3} \le \limsup_{n\to \infty} {X_n^*
     \over (\log n)^3} \le c_2 \; , \qquad \hbox{\rm
     $\p$-a.s.}
     \label{main}
 \end{equation}
\end{theorem}

\medskip

\begin{remark}
 \label{r:main}
{\rm
(i) Theorem \ref{t:main} somehow reminds of Sinai's result (\cite{sinai}) of slow movement of recurrent one-dimensional RWRE, whereas (\ref{yztree}) is a (weaker) analogue of the Kesten--Kozlov--Spitzer characterization (\cite{kesten-kozlov-spitzer}) of sub-diffusive behaviours of {\it transient} one-dimensional RWRE.

(ii) It is interesting to note that tree-valued RWRE possesses both regimes (slow movement and sub-diffusivity) in the recurrent case.

(iii) We mention an important difference between Theorem \ref{t:main} and Sinai's result. If $(Y_n, \; n\ge 0)$ is a recurrent {\it one-dimensional$\,$} RWRE, Sinai's theorem says that ${Y_n \over (\log n)^2}$ converges in distribution (under $\p$) to a non-degenerate limit law, whereas it is known (see \cite{yzsinai}) that
$$
\limsup_{n\to \infty} {Y_n^* \over (\log n)^2} =\infty, \qquad \liminf_{n\to \infty} {Y_n^* \over (\log n)^2} = 0, \qquad \hbox{\rm $\p$-a.s.,}
$$

\noindent where $Y_n^* := \max_{0\le k\le n} |Y_k|$.

(iv) It is not clear to us whether ${X_n^*\over (\log n)^3}$ converges $\p$-almost surely.

(v) We believe that ${|X_n|\over (\log n)^3}$ would converge {\it in distribution} under $\p$.\hfill$\Box$
} 
\end{remark}

\medskip

In Section \ref{s:brw}, we describe the method used to prove Theorem \ref{t:main}. In particular, we introduce an associate branching random walk, and prove an almost sure result for this branching random walk (Theorem \ref{t:V*}) which may be of independent interest. (The two theorems are related to via Proposition \ref{p:V*}.)

The organization of the proof of the theorems is described at the end of Section \ref{s:brw}. We mention that Theorem \ref{t:main} is proved in Section \ref{s:proof}.

Throughout the paper, $c$ (possibly with a subscript) denotes a finite and positive constant; we write $c(\omega)$ instead of $c$ when the value of $c$ depends on the environment $\omega$.

\section{An associated branching random walk}
\label{s:brw}

For any $m\ge 0$, let
$$
\T_m := \left\{ x\in \T: \, |x|=m \right\},
$$

\noindent which stands for the $m$-th generation of the tree. For any $n\ge 0$, let
$$
\tau_n := \inf \left\{ i\ge 1: X_i \in \T_n \right\} = \inf \left\{ i\ge 1: |X_i|= n \right\} ,
$$

\noindent the first hitting time of the walk at level $n$ (whereas $\tau_0$ is the first {\it return} time to the root). We write
$$
\varrho_n := P_\omega \left\{ \tau_n < \tau_0 \right\} .
$$

\noindent In words, $\varrho_n$ denotes the (quenched) probability that the RWRE makes an excursion of height of at least $n$.

An important step in the proof of Theorem \ref{t:main} is the following estimate for $\varrho_n$, in case $\psi'(1)=0$:

\medskip

\begin{theorem}
 \label{t:beta}
 Assume $p={1\over \deg}$ and $\psi'(1)=0$.

 {\rm (i)} There
 exist constants $0<c_3\le c_4<\infty$ such
 that $\P$-almost surely for all large $n$,
 \begin{equation}
     \ee^{-c_4 \, n^{1/3}} \le \varrho_n \le \ee^{-
     c_3 \, n^{1/3}}.
    \label{beta-asymp}
 \end{equation}

 {\rm (ii)} There exist constants $0<c_5\le c_6<
 \infty$ such that for all large $n$,
 \begin{equation}
     \ee^{-c_6 \, n^{1/3}} \le \E (\varrho_n)
     \le \ee^{- c_5 \, n^{1/3}} .
     \label{E(beta)-asymp}
 \end{equation}
\end{theorem}

\medskip

It turns out that $\varrho_n$ is closely related to a branching random walk. But let us first extend the definition of $A(x)$ to all $x\in \T\backslash \{ e\}$.

For any $x\in \T$, let $\{ x_i\}_{1\le i\le \deg}$ denote the set of the children of $x$. In addition of the random variables $A(x)$ ($|x| \ge 2$) defined in (\ref{A}), let $(A(e_i), \; 1\le i\le \deg)$ be a random vector independent of $(\omega(x,y), \, |x|\ge 1, \, y\in \T)$, and distributed as $(A(x_i), \, 1\le i\le \deg)$, for any $x\in \T_m$ with $m\ge 1$. As such, $A(x)$ is well-defined\footnote{The values of $\omega$ at a finite number of edges are of no particular interest. Our choice of $(A(e_i), \; 1\le i\le \deg)$ allows to make unified statements of $A(x)$, $V(x)$, etc., without having to distinguish whether $|x|=1$ or $|x|\ge 2$.} for all $x\in \T \backslash\{ e\}$.

For any $x\in \T\backslash \{ e\}$, the set of vertices on the shortest path relating $e$ and $x$ is denoted by $[\! [ e, \, x] \! ]$; we also set $]\! ] e, \, x] \! ]$ to be $[\! [ e, \, x] \! ] \backslash \{ e\}$.

We now define the process $V=(V(x), \; x\in \T)$ by $V(e) :=0$ and
$$
V(x) := - \sum_{z\in \, ]\!] e,\, x]\!]} \log A(z) , \qquad x\in \T \backslash \{ e\} .
$$

\noindent It is clear that $V$ only depends on the environment $\omega$. In the literature, $V$ is often referred to as a branching random walk, see for example Biggins and Kyprianou~\cite{biggins-kyprianou97}.

We first state the main result of the section. Let
\begin{equation}
    \overline{V} (x) := \max_{z\in \, ]\!] e, \,
    x]\!]} V(z) ,
    \label{Vbar}
\end{equation}

\noindent which stands for the maximum of $V$ over the path $]\!] e, \, x]\!]$.

\medskip

\begin{theorem}
 \label{t:V*}
 If $p={1\over \deg}$ and if $\psi'(1)\ge 0$, then
 there exist constants $0<c_7 \le c_8<\infty$ such
 that
 \begin{equation}
     c_7 \le \liminf_{n\to \infty} {1\over n^{1/3}}
     \min_{x\in \T_n} \overline{V} (x) \le
     \limsup_{n\to \infty} {1\over n^{1/3}}\min_{x
     \in \T_n} \overline{V} (x) \le c_8 \, , \qquad
     \hbox {\rm $\P$-a.s.}
     \label{V*}
 \end{equation}
\end{theorem}

\medskip

\begin{remark}
 \label{r:V*}
{\rm
(i) We cannot replace $\min_{x\in \T_n} \overline{V} (x)$ by $\min_{x\in \T_n} V(x)$ in Theorem \ref{t:V*}; in fact, it is proved by McDiarmid~\cite{mcdiarmid} that there exists a constant $c_9$ such that $\P$-almost surely for all large $n$, we have $\min_{x\in \T_n} V(x) \le c_9\log n$.

(ii) If ($p={1\over \deg}$ and) $\psi'(1)< 0$, it is well-known (Hammersley~\cite{hammersley}, Kingman~\cite{kingman}, Biggins~\cite{biggins}) that ${1\over n}\min_{x\in \T_n} V(x)$ converges $\P$-almost surely to a (strictly) positive constant whose value is known; thus $\min_{x\in \T_n} \overline{V}(x)$ grows linearly in this case.

(iii) Only the case $\psi'(1)=0$ needs to be proved. Indeed, if ($p={1\over \deg}$ and) $\psi'(1)>0$, then there exists a unique $0<\theta<1$ such that
\begin{equation}
    \psi'(\theta) =0, \qquad \E(A^\theta) = {1\over
    \deg} .
    \label{gamma}
\end{equation}

\noindent We define $\widetilde{A} := A^\theta$, $\widetilde{p} := \inf_{t\in [0,1]} \E(\widetilde{A}^t)$ and $\widetilde{\psi} (t) := \log \E(\widetilde{A}^t)$, $t\ge 0$. Clearly, we have
$$
\widetilde{p}={1\over \deg}, \qquad \widetilde{\psi}'(1)=0.
$$

\noindent Let $\widetilde{V} (x) := -\sum_{z\in \, ]\!] e, \, x]\!]} \log \widetilde{A}(z)$. Then $V(x) = {1\over \theta}\, \widetilde{V} (x)$, which leads us to the case $\psi'(1)=0$.\hfill$\Box$
} 
\end{remark}

\medskip

Here is the promised relation between $\varrho_n$ and $V$, for recurrent RWRE on $\T$.

\medskip

\begin{proposition}
 \label{p:V*}
 If $(X_n)$ is recurrent, there exists a constant
 $c_{10}>0$ such that for any $n \ge 1$,
 \begin{equation}
     \varrho_n \ge {c_{10}\over n} \exp \left( -
     \min_{x \in \T_n} \overline{V}(x) \right) .
    \label{beta>V}
 \end{equation}
\end{proposition}

\medskip

\noindent {\it Proof of Proposition \ref{p:V*}.} For any $x\in \T$, let
\begin{equation}
    T(x) := \inf \left\{ i\ge 0: X_i =x \right\},
    \label{T(x)}
\end{equation}

\noindent which is the first hitting time of the walk at vertex $x$. By definition, $\tau_n = \min_{x \in \T_n} T(x)$, for $n\ge 1$. Therefore,
\begin{equation}
    \varrho_n \ge \max_{x\in \T_n} P_\omega \left\{
    T(x) < \tau_0 \right\} .
    \label{beta>T}
\end{equation}

We now compute the (quenched) probability $P_\omega \{ T(x) < \tau_0 \}$. We fix $x\in \T_n$, and define a random sequence $(\sigma_j)_{j\ge 0}$ by $\sigma_0:= 0$ and
$$
\sigma_j:= \inf \left\{ k> \sigma_{j-1}: X_k \in [\! [e,\, x] \! ] \backslash \{ X_{\sigma_{j-1}} \} \right\}, \qquad j\ge 1.
$$

\noindent (Of course, the sequence depends on $x$.) Let
\begin{equation}
    Z_k := X_{\sigma_k}, \qquad k\ge 0.
    \label{Xhat}
\end{equation}

\noindent In words, $Z=(Z_k, \; k\ge 0)$ is the restriction of $X$ on the path $[\! [e,\, x] \! ]$; i.e., it is almost the original walk, except that we remove excursions away from $[\! [e,\, x] \! ]$. Clearly, $Z$ is a one-dimensional RWRE with (writing $[\! [e,\, x ] \! ] = \{ e=: x^{(0)}, x^{(1)}, \cdots , x^{(n)} := x\}$)
\begin{eqnarray*}
    P_\omega \left\{ Z_{k+1}= x^{(i+1)}
    \, \Big| \, Z_k= x^{(i)} \right\}
 &=& {A(x^{(i+1)}) \over 1+ A(x^{(i+1)})},
    \\
    P_\omega \left\{ Z_{k+1}= x^{(i-1)}
    \, \Big| \, Z_k= x^{(i)} \right\}
 &=& {1\over 1+ A(x^{(i+1)})},
\end{eqnarray*}

\noindent for all $1\le i\le n-1$. We observe that
\begin{eqnarray*}
    P_\omega \{ T(x) < \tau_0 \}
 &=& \omega(e, x^{(1)}) \, P_\omega \left\{
    \hbox{$Z$ hits $x^{(n)}$ before
    hitting $e$} \, \Big| \, Z_0 = x^{(1)}
    \right\}
    \\
 &=& \omega(e, x^{(1)}) \, {\ee^{V(x^{(1)})} \over
    \sum_{z\in \, ]\! ] e,\, x] \! ]} \ee^{V(z)} } ,
\end{eqnarray*}

\noindent the second identity following from a general formula (Zeitouni~\cite{zeitouni}, formula (2.1.4)) for the exit problem of one-dimensional RWRE. By ellipticity condition, there exists a constant $c_{11}>0$ such that $\omega(e, x^{(1)}) \ee^{V(x^{(1)})} \ge c_{11}$. Plugging this estimate into (\ref{beta>T}) yields
$$
\varrho_n \ge \max_{x\in \T_n} {c_{11} \over \sum_{y\in \, ]\! ] e, \, x]\! ]} \ee^{V(y)}} ,
$$

\noindent completing the proof of Proposition \ref{p:V*}.\hfill$\Box$

\bigskip

The proof of the theorems is organized as follows.
\begin{itemize}

\item{Section \ref{s:V-ub}: Theorem \ref{t:V*}, upper bound}.

\item{Section \ref{s:beta}: Theorem \ref{t:beta} (by means of the upper bound in Theorem \ref{t:V*}; this is the technical part of the paper)}.

\item{Section \ref{s:V-lb}: Theorem \ref{t:V*}, lower bound (by means of the upper bound in Theorem \ref{t:beta})}.

\item{Section \ref{s:proof}: Theorem \ref{t:main}}.

\end{itemize}

\section{Proof of Theorem \ref{t:V*}: upper bound}
\label{s:V-ub}

Throughout the section, we assume $p={1\over \deg}$ and $\psi'(1)=0$.

Let
\begin{equation}
    B(x) := \prod_{y\in \, ] \! ]e, \, x ]\! ]} A
    (y), \qquad x\in \T \backslash\{ e\}.
    \label{B}
\end{equation}

\noindent We start by recalling a change-of-probability formula from Biggins and Kyprianou \cite{biggins-kyprianou97}. See also Durrett and Liggett~\cite{durrett-liggett}, and Bingham and Doney~\cite{bingham-doney}.

\medskip

\begin{fact} {\bf (Biggins and Kyprianou \cite{biggins-kyprianou97}).}
 \label{f:biggins-kyprianou}
 For any $n\ge 1$ and any positive measurable
 function $G$,
 \begin{equation}
    \sum_{x\in \T_n} \E \left[ B(x)\, G\left( B(z),
    \; z\in \, ]\! ] e, \, x] \!] \right) \right] =
    \E \left[ G\left( \ee^{S_i} , \; 1\le i\le n
    \right) \right] ,
    \label{biggins}
\end{equation}

\noindent where $S_n$ is the sum of $n$ i.i.d.\ centered random variables whose common distribution is determined by
$$
\E \left[ g\left( S_1\right) \right] = \deg \, \E \left[ A \, g\left( \log A \right) \right] ,
$$

\noindent for any positive measurable function $g$.
\end{fact}

\medskip

The fact that $S_1$ is centered is a consequence of the assumption $\psi'(1)=0$. We note that in (\ref{biggins}), the value of $\E [ B(x)\, G ( B(z), \; z\in \, ]\! ] e, \, x] \!] ) ]$ is the same for all $x\in \T_n$.

We have now all the ingredients of the proof of the upper bound in Theorem \ref{t:V*}.

\bigskip

\noindent {\it Proof of Theorem \ref{t:V*}: upper bound.} By Remark \ref{r:V*}, only the case $\psi'(1)=0$ needs to be treated. We assume in the rest of the section ($p={1\over \deg}$ and) $\psi'(1)=0$. The proof borrows some ideas of Bramson~\cite{bramson} concerning branching Brownian motions. Let
$$
E_m := \left\{ x\in \T_m: \; \max_{z\in \, ]\!] e, \, x]\!]} |V(z)| \le m^{1/3} \right\}.
$$

We first estimate $\E [\# E_m]$:
$$
\E [\# E_m] = \sum_{x\in \T_m} \P \left\{ \, \max_{z\in \, ]\!] e, \, x]\!]} |V(z)| \le m^{1/3} \right\} .
$$

\noindent By assumption, for any given $x\in \T_m$, $(V(z), \; z\in \, ]\!] e, \, x]\!])$ is the set of the first $m$ partial sums of i.i.d.\ random variables whose common distribution is $A$. By (\ref{biggins}), this leads to:
$$
\E [\# E_m] = \E \left( \ee^{-S_m} \, {\bf 1}_{ \{ \max_{1\le i\le m} |S_i| \le m^{1/3} \} } \right) \ge \P \left\{ \max_{1\le i\le m} |S_i| \le m^{1/3}, \; S_m \le 0 \right\} .
$$

\noindent The probability on the right-hand side is a ``small deviation" probability, with an unimportant condition upon the terminal value. By a general result of Mogul'skii~\cite{mogulskii}, we have, for all sufficiently large $m$ (say $m\ge m_0$),
$$
\E [\# E_m] \ge \exp \left( - c_{12} \, m^{1/3} \right) .
$$

We now estimate the second moment of $\# E_m$. For any pair of vertices $x$ and $y$, we write $x<y$ if $x$ is an ancestor of $y$, and $x\le y$ if $x$ is either $y$ itself or an ancestor of $y$. Then
\begin{eqnarray*}
 && \E [(\# E_m)^2] - \E [\# E_m]
    \\
 &=& \sum_{u, v\in \T_m, \; u\not= v} \P \left\{ u
    \in E_m, \; v\in E_m \right\}
    \\
 &=& \sum_{j=0}^{m-1}\; \sum_{z\in \T_j}\;
    \sum_{x\in \T_{j+1}: \; z<x} \; \sum_{y\in
    \T_{j+1}\backslash\{x\}: \; z<y} \;
    \sum_{u\in \T_m: \; x\le u} \; \sum_{v\in \T_m:
    \; y\le v}\P \left\{ u \in E_m, \; v\in E_m
    \right\}.
\end{eqnarray*}

\noindent In words, $z$ is the youngest common ancestor of $u$ and $v$, while $x$ and $y$ are distinct children of $z$ at generation $j+1$. If $j=m-1$, we have $x=u$ and $y=v$, otherwise $x$ is an ancestor of $u$, and $y$ of $v$.

Fix $z\in \T_j$, and let $x$ and $y$ be a pair of distinct children of $z$. Let $u\in \T_m$ and $v\in \T_m$ be such that $x\le u$ and $y\le v$. Then
\begin{eqnarray*}
 &&\P \left\{ u \in E_m, \; v\in E_m\right\}
    \\
 &\le& \P \left\{ \max_{r\in \, ]\!] e, \, z]\!]}
    |V(r)| \le m^{1/3} \right\} \times \left( \P
    \left\{ \max_{r\in \, ]\!] z, \, x]\!]} |V(r) -
    V(z)| \le 2m^{1/3} \right\} \right) ^2.
\end{eqnarray*}

\noindent We have, by (\ref{biggins}),
$$
\P \left\{ \max_{r\in \, ]\!] e, \, z]\!]} |V(r)| \le m^{1/3} \right\} = \deg^{-j}\, \E \left[ \ee^{-S_j} \, {\bf 1}_{ \{ \max_{1\le i\le j} |S_i| \le
m^{1/3} \} } \right] \le \deg^{-j} \ee^{m^{1/3}} ,
$$

\noindent and similarly, $\P \{ \max_{r\in \, ]\!] z, \, x]\!]} |V(r) - V(z)| \le 2m^{1/3} \} \le \deg^{-(m-j)}\ee^{2m^{1/3}}$. Therefore,
\begin{eqnarray*}
 && \E [(\# E_m)^2] - \E [\# E_m]
    \\
 &\le& \sum_{j=0}^{m-1}\; \sum_{z\in \T_j}\;
    \sum_{x\in \T_{j+1}: \; z<x} \; \sum_{y\in
    \T_{j+1}\backslash\{x\}: \; z<y} \;
    \sum_{u\in \T_m: \; x\le u} \; \sum_{v\in \T_m:
    \; y\le v} \deg^{j-2m} \ee^{5m^{1/3}}
    \\
 &=& \sum_{j=0}^{m-1}\; \sum_{z\in \T_j} \deg
    (\deg-1) \deg^{m-j-1} \deg^{m-j-1} \deg^{j-2m}
    \ee^{5m^{1/3}}
    \\
 &=& {\deg-1\over \deg} \, m\, \ee^{5m^{1/3}} .
\end{eqnarray*}

\noindent Recall that $\E [\# E_m] \ge \exp ( - c_{12} \, m^{1/3})$ for $m\ge m_0$. Therefore, for $m\ge m_0$,
$$
{\E [(\# E_m)^2] \over \{ \E [\# E_m]\}^2} \le {\deg-1\over \deg} \, m\, \ee^{(5+2c_{12}) m^{1/3}} + \ee^{c_{12}\, m^{1/3}} \le \ee^{c_{13}\, m^{1/3}} .
$$

\noindent By the Cauchy--Schwarz inequality, for $m\ge m_0$,
$$
\P \left\{ E_m \not= \emptyset \right\} = \P \left\{ \# E_m >0 \right\} \ge {\{ \E [\# E_m]\}^2 \over \E [(\# E_m)^2] } \ge \ee^{-c_{13}\, m^{1/3}}.
$$

\noindent A fortiori, for $m\ge m_0$,
$$
\P \left\{ \exists x\in \T_m, \; \overline{V}(x) \le m^{1/3} \right\} \ge \ee^{-c_{13}\, m^{1/3}},
$$

\noindent which implies
$$
\P \left\{ \min_{x\in \T_m} \overline{V}(x) > m^{1/3} \right\} \le 1- \ee^{-c_{13}\, m^{1/3}} \le \exp\left( - \ee^{-c_{13}\, m^{1/3}} \right) .
$$

\noindent Let $n>m$. By the ellipticity condition stated in the Introduction, there exists a constant $c_{14}>0$ such that $\max_{z\in \, ] \! ] e, \, y] \! ]} V(z) \le c_{14}\, (n-m)$ for any $y\in \T_{n-m}$. Accordingly, for $m\ge m_0$,
\begin{eqnarray*}
 &&\P \left\{ \min_{x\in \T_n} \overline{V}(x) >
    m^{1/3} + c_{14}\, (n-m) \right\}
    \\
 &\le& \P \left\{ \min_{y\in T_{n-m}} \min_{x\in
    \T_n: \, y<x} \max_{r\in \, ] \! ] y, \, x]\! ]}
    [ V(r)- V(y)] > m^{1/3} \right\}
    \\
 &=& \left( \P \left\{ \min_{s\in \T_m}
    \overline{V}(s) > m^{1/3} \right\}
    \right)^{\! \deg^{n-m}}
    \\
 &\le& \exp\left( - \deg^{n-m} \ee^{-c_{13}\,
    m^{1/3}} \right) .
\end{eqnarray*}

\noindent We now choose $m= m(n) := n- \lfloor c_{15} \, n^{1/3} \rfloor$, where the constant $c_{15}$ is sufficiently large such that $\sum_n \exp ( - \deg^{n-m} \ee^{-c_{13}\, m^{1/3}} )<\infty$. Then, by the Borel--Cantelli lemma,
$$
\limsup_{n\to \infty}\, {1\over n^{1/3}} \min_{x\in \T_n} \overline{V}(x) \le 1+ c_{14}c_{15}, \qquad \hbox{\rm $\P$-a.s.,}
$$

\noindent yielding the desired upper bound in Theorem \ref{t:V*}.\hfill$\Box$

\section{Proof of Theorem \ref{t:beta}}
\label{s:beta}

Throughout the section, we assume $p={1\over \deg}$ and $\psi'(1)=0$.

\bigskip

\noindent {\it Proof of Theorem \ref{t:beta}: lower bound.} The estimate $\varrho_n \ge \ee^{-c_4 \, n^{1/3}}$ ($\P$-almost surely for all large $n$) follows immediately from the upper bound in Theorem \ref{t:V*} (proved in Section \ref{s:V-ub}) by means of Proposition \ref{p:V*}, with any constant $c_4> c_8$. By Fatou's lemma, we have $\liminf_{n\to \infty} \ee^{c_4 \, n^{1/3}} \E(\varrho_n) \ge 1$.\hfill$\Box$

\bigskip

We now introduce the important ``additive martingale" $M_n$; in particular, the lower tail behaviour of $M_n$ is studied in Lemma \ref{l:ltail}, by means of another martingale called ``multiplicative martingale". The upper bound in Theorem \ref{t:beta} will then be proved based on the asymptotics of $M_n$ and on the just proved lower bound.

Let $B(x) := \prod_{y\in \, ] \! ]e, \, x ]\! ]} A(y)$ (for $x\in \T \backslash\{ e\}$) as in (\ref{B}), and let
\begin{equation}
    M_n := \sum_{x\in \T_n} B(x), \qquad n\ge 1.
    \label{Mn}
\end{equation}

\noindent When $\E(A) = {1\over \deg}$ (which is the case if $p={1\over \deg}$ and $\psi'(1)=0$), the process $(M_n, \; n\ge 1)$ is a martingale, and is referred to as an associated ``{\it additive martingale}$\,$".

It is more convenient to study the behaviour of $M_n$ by means of another martingale. It is known (see Liu~\cite{liu00}) that under assumptions $p={1\over \deg}$ and $\psi'(1)=0$, there is a unique non-trivial function $\varphi^*: \r_+ \to (0, \, 1]$ such that
\begin{equation}
    \varphi^*(t) = \E \left\{ \prod_{i=1}^\deg
    \varphi^* (t A(e_i)) \right\}, \qquad t\ge 0.
    \label{cascade}
\end{equation}

\noindent (By non-trivial, we mean that $\varphi^*$ is not identically 1.) Let
$$
M_n^* := \prod_{x\in \T_n} \varphi^*(B(x)), \qquad n\ge 1.
$$

\noindent The process $(M_n^*, \; n\ge 1)$ is also a martingale (Liu~\cite{liu00}). Following Neveu~\cite{neveu}, we call $M_n^*$ an associated ``{\it multiplicative martingale}$\,$".

Since the martingale $M_n^*$ takes values in $(0,\, 1]$, it converges almost surely (when $n\to \infty$) to, say, $M_\infty^*$, and $\E(M_\infty^*)=1$. It is proved by Liu~\cite{liu00} that $\E \{ (M_\infty^*)^t\} = \varphi^* (t)$ for any $t\ge 0$.

Recall that for some $0<\alpha <1$,
\begin{eqnarray}
    \log \left( {1\over \varphi^*(t)} \right)
 &\sim& t \log \left( {1\over t}\right) , \qquad t
    \to 0,
    \label{phi(0)}
    \\
    \log \left( {1\over \varphi^* (s)} \right)
 &\ge& c_{16} \, s^\alpha, \qquad s\ge 1;
    \label{phi(infty)1}
\end{eqnarray}

\noindent see Liu (\cite{liu00}, Theorem 2.5) for (\ref{phi(0)}), and Liu (\cite{liu01}, Theorem 2.5) for (\ref{phi(infty)1}).

\medskip

\begin{lemma}
 \label{l:ltail}
 Assume $p={1\over \deg}$ and $\psi'(1)=0$. For any
 $\chi > 1/2$, there exists $\delta>0$ such that
 for all sufficiently large $n$,
 \begin{equation}
     \P \left\{ M_n < n^{-\chi} \right\} \le \exp
     \left( - n^\delta \right) .
     \label{M-ltail}
 \end{equation}

\end{lemma}

\medskip

\noindent {\it Proof of Lemma \ref{l:ltail}.} Let $K>0$ be such that $\P\{ M_\infty^*> \ee^{- K} \} >0$. Then $\varphi^* (t) = \E \{ (M_\infty^*)^t\} \ge \P\{ M_\infty^*> \ee^{- K} \} \, \ee^{-K \, t}$ for all $t>0$. Thus, there exists $c_{17}>0$ such that for all $t\ge 1$, $\varphi^* (t) \ge \ee^{-c_{17} \, t}$.

Let $\varepsilon >0$. By (\ref{phi(0)}) and (\ref{phi(infty)1}), there exists a constant $c_{18}$ such that
$$
\log \left( {1\over M_n^*}\right) = \sum_{x\in \T_n} \log \left( {1\over \varphi^*(B(x))}\right) \le c_{18} \left( J_{1,n} + J_{2,n} + J_{3,n} \right) ,
$$

\noindent where
\begin{eqnarray*}
    J_{1,n}
 &:=& \sum_{x\in \T_n} B(x) \left( \, \log {1
    \over B (x)} \right) \, {\bf 1}_{ \{ B(x) <
    \exp(-n^{(1/2)+\varepsilon}) \} } ,
    \\
    J_{2,n}
 &:=& \sum_{x\in \T_n} B(x) \left( \, \log {\ee
    \over B (x)} \right) \, {\bf 1}_{
    \{ \exp(-n^{(1/2)+\varepsilon}) \le B (x) \le 1
    \} } ,
    \\
    J_{3,n}
 &:=& \sum_{x\in \T_n} B(x) \, {\bf 1}_{
    \{ B(x) > 1\} } .
\end{eqnarray*}

\noindent Clearly, $J_{3,n} \le \sum_{x\in \T_n} B(x) = M_n$, whereas $J_{2,n} \le (n^{(1/2)+\varepsilon}+1) M_n$. Hence, $J_{2,n} + J_{3,n} \le (n^{(1/2)+\varepsilon}+2) M_n \le 2n^{(1/2)+\varepsilon} M_n$ (for $n\ge 4$). Accordingly, for $n\ge 4$,
\begin{equation}
    n^{(1/2)+\varepsilon} M_n \ge {1\over 2c_{18}}
    \log \left( {1\over M_n^*}\right) - {1\over 2}
    J_{1,n}.
    \label{Mn*}
\end{equation}

We now estimate the tail probability of $M_n^*$. Let $\lambda\ge 1$ and $z>0$. By Chebyshev's inequality,
$$
\P\left\{ \log \left( {1\over M_n^*}\right) < z\right\} \le \ee^{\lambda z}\, \E \left\{ (M_n^*)^\lambda \right\}.
$$

\noindent Since $M_n^*$ is a bounded martingale, $\E \{ (M_n^*)^\lambda \} \le \E \{ (M_\infty^*)^\lambda \} = \varphi^* (\lambda)$. Therefore,
$$
\P\left\{ \log \left( {1\over M_n^*}\right) < z\right\} \le \ee^{\lambda z}\varphi^* (\lambda).
$$

\noindent Choosing $z:= 4c_{18}\, n^{-\varepsilon}$ and $\lambda:= n^\varepsilon$, it follows from (\ref{phi(infty)1}) that
$$
\P\left\{ \log \left( {1\over M_n^*}\right) <
4c_{18}\, n^{-\varepsilon} \right\} \le \exp\left(
4c_{18} - c_{16} \; n^{\varepsilon \alpha} \right) .
$$

\noindent Plugging this into (\ref{Mn*}) yields that for $n\ge 4$,
\begin{equation}
    \P \left\{ n^{(1/2)+\varepsilon} M_n + {1\over
    2} J_{1,n} < 2n^{-\varepsilon} \right\} \le
    \exp\left( 4c_{18} - c_{16} \; n^{\varepsilon
    \alpha} \right) .
    \label{Mn>>}
\end{equation}

We note that $J_{1,n} \ge 0$. By (\ref{biggins}),
$$
\E(J_{1,n}) = \E \left\{ (-S_n) \, {\bf 1}_{ \{ S_n <- n^{(1/2) + \varepsilon} \} } \right\} .
$$

\noindent Recall that $S_n$ is the sum of $n$ i.i.d.\ bounded centered random variables. It follows that for all sufficiently large $n$,
$$
\E(J_{1,n}) \le \exp \left( - c_{19} \, n^{2\varepsilon} \right) .
$$

\noindent By (\ref{Mn>>}) and Chebyshev's inequality,
\begin{eqnarray*}
    \P \left\{ n^{(1/2)+\varepsilon} M_n < n^{-
    \varepsilon} \right\}
 &\le & \P \left\{ n^{(1/2)+\varepsilon} M_n + {1
    \over 2} J_{1,n} < 2n^{-\varepsilon} \right\} +
    \P \left\{ J_{1,n} \ge 2n^{-\varepsilon}
    \right\}
    \\
 &\le & \exp\left( 4c_{18} - c_{16} \;
    n^{\varepsilon \alpha} \right) +
    {n^\varepsilon \over 2} \exp \left( - c_{19} \,
    n^{2\varepsilon} \right) ,
\end{eqnarray*}

\noindent from which (\ref{M-ltail}) follows.\hfill$\Box$

\bigskip

We have now all the ingredients for the proof of the upper bound in Theorem \ref{t:beta}.

\bigskip

\noindent {\it Proof of Theorem \ref{t:beta}: upper bound.} We only need to prove the upper bound in (\ref{E(beta)-asymp}), namely, there exists $c_5$ such that for all large $n$,
\begin{equation}
    \E (\varrho_n) \le \ee^{ - c_5 \, n^{1/3}} .
    \label{E(beta)<<}
\end{equation}

\noindent If (\ref{E(beta)<<}) holds, then the upper bound in (\ref{beta-asymp}) follows by an application of Chebyshev's inequality and the Borel--Cantelli lemma.

It remains to prove (\ref{E(beta)<<}). For any $x\in \T \backslash \{ e\}$, we define
$$
\beta_n(x) := P_\omega \left\{ \hbox{starting from $x$, the RWRE hits $\T_n$ before hitting ${\buildrel
\leftarrow \over x}$} \right\} ,
$$

\noindent where, as before, ${\buildrel \leftarrow \over x}$ is the parent of $x$. In the notation of (\ref{T(x)}),
$$
\beta_n(x) = P_\omega \{ T_n < T({\buildrel \leftarrow \over x}) \, | \, X_0=x \} ,
$$

\noindent where $T_n := \min_{x\in \T_n} T(x)$. Clearly, $\beta_n(x) =1$ if $x\in \T_n$.

Recall that for any $x\in \T$, $\{ x_i\}_{1\le i\le \deg}$ is the set of the children of $x$. By the Markov property, if $1\le |x| \le n-1$, then
$$
\beta_n(x) = \sum_{i=1}^\deg \omega(x, x_i) P_\omega \{ T_n < T({\buildrel \leftarrow \over x}) \, | \, X_0=x_i \}.
$$

\noindent Consider the event $\{ T_n < T({\buildrel \leftarrow \over x})\}$ when the walk starts from $x_i$. There are two possible situations: (i) either $T_n< T(x)$ (which happens with probability $\beta_n(x_i)$, by definition); (ii) or $T_n>T(x)$ and after hitting $x$ for the first time, the walk hits $\T_n$ before hitting ${\buildrel \leftarrow \over x}$. By the strong Markov property, $P_\omega \{ T_n < T({\buildrel \leftarrow \over x}) \, | \, X_0=x_i \} = \beta_n(x_i) + [1-\beta_n(x_i)] \beta_n(x)$. Therefore,
\begin{eqnarray*}
    \beta_n(x)
 &=& \sum_{i=1}^\deg \omega(x, x_i) \beta_n(x_i) +
    \beta_n(x) \sum_{i=1}^\deg \omega(x, x_i)
    [1-\beta_n(x_i)]
    \\
 &=& \sum_{i=1}^\deg \omega(x, x_i) \beta_n(x_i) +
    \beta_n(x) [1- \omega(x, {\buildrel \leftarrow
    \over x})] - \beta_n(x) \sum_{i=1}^\deg \omega
    (x, x_i) \beta_n(x_i) ,
\end{eqnarray*}

\noindent from which it follows that
\begin{equation}
    \beta_n(x) = {\sum_{i=1}^\deg A(x_i) \beta_n
    (x_i) \over 1+ \sum_{i=1}^\deg A(x_i) \beta_n
    (x_i) }, \qquad 1\le |x| \le n-1.
    \label{beta-recurrence}
\end{equation}

\noindent Together with condition $\beta_n(x) =1$ (for $x\in \T_n$), these equations determine the value of $\beta_n(x)$ for all $x\in \T$ such that $1\le |x| \le n$.

We introduce the random variable
\begin{equation}
    \beta_n(e) := {\sum_{i=1}^\deg A(e_i) \beta_n
    (e_i) \over 1+ \sum_{i=1}^\deg A(e_i) \beta_n
    (e_i)} .
    \label{beta_n(e)}
\end{equation}

\noindent The value of $\beta_n(e)$ for given $\omega$ is of no importance, but the distribution of $\beta_n(e)$, which is identical to that of $\beta_{n+1}(e_1)$, plays a certain role in several places of the proof. For example, for $1\le |x| < n$, the random variables $\beta_n(x)$ and $\beta_{n-|x|}(e)$ have the same distribution; in particular, $\E[\beta_n(x)] = \E[\beta_{n-|x|}(e)]$. In the rest of the section, we make frequent use of this property without further mention. We also make the trivial observation that for $1\le |x| < n$, $\beta_n(x)$ depends only on those $A(y)$ such that $|x|+1 \le |y| \le n$ and that $x$ is an ancestor of $y$.

Recall that $\varrho_n = P_\omega \{ \tau_n < \tau_0 \}$. Therefore,
\begin{equation}
    \varrho_n = \sum_{i=1}^\deg \omega(e, \, e_i)
    \beta_n(e_i).
    \label{beta-rho}
\end{equation}

\noindent In particular,
\begin{equation}
    \E(\varrho_n) = \E[\beta_n(e_i)] = \E[
    \beta_{n-1}(e)] , \qquad \forall 1\le i\le \deg.
    \label{E(beta)=E(beta)}
\end{equation}

Let $a_j := \E(\varrho_{j^3 +1}) = \E[\beta_{j^3}(e)]$, $j=0$, $1$, $2$, $\cdots$, $\lfloor n^{1/3}\rfloor$. Clearly, $a_0=1$, and $j\mapsto a_j$ is non-increasing for $0\le j\le \lfloor n^{1/3}\rfloor$. We look for an upper bound for $a_{\lfloor n^{1/3}\rfloor}$.

Let $m> \Delta \ge 1$ be integers. Let $1\le i\le \deg$, and let $(e_{ij}, \; 1\le j\le \deg)$ be the set of children of $e_i$. By (\ref{beta-recurrence}), we have
$$
\beta_m(e_i) \le \sum_{j=1}^\deg A(e_{ij}) \beta_m (e_{ij}).
$$

\noindent Iterating the same argument, we arrive at:
$$
\beta_m(e_i) \le \sum_{y\in \T_\Delta: \, y<e_i} \left( \prod_{z: \, e_i<z, \, z\le y} A(z) \right) \beta_m(y) = \sum_{y\in \T_\Delta: \, y<e_i} {B(y)\over A(e_i)} \beta_m(y).
$$

\noindent By (\ref{beta_n(e)}), this yields
$$
\beta_m(e) \le {\sum_{i=1}^\deg \sum_{y\in \T_\Delta: \, y<e_i} B(y) \beta_m(y) \over 1+ \sum_{i=1}^\deg \sum_{y\in \T_\Delta: \, y<e_i} B(y) \beta_m(y)} ={\sum_{y\in \T_\Delta} B(y) \beta_m(y) \over 1+ \sum_{y\in \T_\Delta} B(y) \beta_m(y)}.
$$

Fix $n$ and $0\le j\le \lfloor n^{1/3}\rfloor -1$. Let
$$
\Delta = \Delta(j) := (j+1)^3-j^3 = 3j^2+3j+1.
$$

\noindent Then
$$
a_{j+1} = \E[\beta_{(j+1)^3}(e)] \le \E \left( {\sum_{y\in \T_{\Delta}} B(y) \beta_{(j+1)^3}(y) \over 1+ \sum_{y\in \T_{\Delta}} B(y) \beta_{(j+1)^3}(y)} \right) .
$$

\noindent We note that $(\beta_{(j+1)^3}(y), \, y\in \T_{\Delta})$ is a collection of i.i.d.\ random variables distributed as $\beta_{j^3}(e)$, and is independent of $(B(y), \, y\in \T_{\Delta})$.

Let $(\xi(x), \; x\in \T)$ be i.i.d.\ random variables distributed as $\beta_{j^3}(e)$, independent of all other random variables and processes. Let
$$
N_m := \sum_{x\in \T_m} B(x) \xi(x), \qquad m\ge 1.
$$

\noindent The last inequality can be written as
\begin{equation}
    a_{j+1} \le \E \left( {N_{\Delta}\over 1+
    N_{\Delta} }\right) .
    \label{aj}
\end{equation}

\noindent By definition,
\begin{equation}
    \E \left( {N_{\Delta} \over 1+ N_{\Delta}}
    \right) = \sum_{x\in \T_{\Delta}} \E \left(
    {B(x) \xi(x) \over 1+ N_{\Delta}} \right) =
    \sum_{x\in \T_{\Delta}} \E \left\{ B(x) \xi(x)
    \ee^{-Y N_{\Delta}} \right\},
    \label{toto1}
\end{equation}

\noindent where $Y$ is an exponential random variable of parameter 1, independent of everything else.

\noindent Let us fix $x\in \T_{\Delta}$, and estimate $\E\{ B(x) \xi(x) \ee^{-Y N_{\Delta}} \}$. Since $N_m=\sum_{x\in \T_m} B(x) \xi(x)$ (for any $m\ge 1$), we have
$$
N_\Delta \ge B({\buildrel \leftarrow \over x} ) A(y) \xi(y) ,
$$

\noindent for any $y\in \T_\Delta \backslash \{ x\}$ such that ${\buildrel \leftarrow \over y} = {\buildrel \leftarrow \over x}$. Note that by ellipticity condition, $A(y) \ge c>0$ for some constant $c$. Accordingly,
\begin{eqnarray*}
    \E \left\{ B(x) \xi(x) \ee^{-Y N_{\Delta}}
    \right\}
 &\le& \E \left\{ B(x) \xi(x) \ee^{-cY B(
    {\buildrel \leftarrow \over x} ) \xi(y) }
    \right\}
    \\
 &=& \E \left\{ \xi(x) \right\} \E \left\{ B(x)
    \ee^{-cY B({\buildrel \leftarrow \over x} )
    \xi(y) }\right\}.
\end{eqnarray*}

\noindent Recall that $\e \{ \xi(x) \} = \E \{ \beta_{j^3}(e) \} = a_j$ and that $\xi(y)$ is distributed as $\beta_{j^3}(e)$, independent of $(B(x), Y,B({\buildrel \leftarrow \over x} ) )$. At this stage, it is convenient to recall the following inequality (see \cite{yztree} for an elementary proof): if $\E(A)= {1\over \deg}$ (which is guaranteed by the assumption $p={1\over \deg}$ and $\psi'(1)=0$), then
$$
\E \left\{ \exp \left( -t \, {\beta_k (e) \over \E [\beta_k(e)]} \right) \right\} \le \E \left\{ \ee^{- t M_k} \right\}, \qquad \forall k\ge 1, \; \; \forall t \ge 0,
$$

\noindent where $M_k$ is defined in (\ref{Mn}). As a consequence,
$$
\E \left\{ B(x) \xi(x) \ee^{-Y N_{\Delta}} \right\} \le a_j \, \E \left\{ B(x) \ee^{-cY B({\buildrel \leftarrow \over x} ) \, a_j \widetilde{M}_{j^3} }\right\} ,
$$

\noindent where $\widetilde{M}_{j^3}$ is distributed as $M_{j^3}$, and is independent of everything else.
Since $A(x)= {B(x) \over B ( {\buildrel \leftarrow \over x})}$ is independent of $B ({\buildrel \leftarrow \over x})$ (and $Y$ and $\widetilde{M}_{j^3}$), with $\E \{ A(x) \} = {1\over \deg}$, this yields
$$
\E \left\{ B(x) \xi(x) \ee^{-Y N_{\Delta}} \right\} \le {a_j \over \deg} \E \left\{ B({\buildrel \leftarrow \over x}) \ee^{- c\, a_j\, Y B ( {\buildrel \leftarrow \over x}) \widetilde{M}_{j^3} } \right\} .
$$

\noindent Plugging this into (\ref{toto1}), we see that
\begin{eqnarray*}
    \E \left( {N_{\Delta} \over 1+ N_{\Delta}}
    \right)
 &\le& a_j \sum_{u\in \T_{\Delta-1}} \E \left\{ B(u)
    \ee^{- c\, a_j \, Y B (u)\widetilde{M}_{j^3} }
    \right\}
    \\
 &=& a_j \, \E \left\{ \exp\left( - c\, a_j \, Y\,
    \ee^{S_{\Delta-1}} \widetilde{M}_{j^3} \right)
    \right\} ,
\end{eqnarray*}

\noindent the last identity being a consequence of (\ref{biggins}), the random variables $Y$, $S_{\Delta-1}$ and $\widetilde{M}_{j^3}$ being independent. By (\ref{aj}), $a_{j+1} \le \E ( {N_{\Delta}\over 1+ N_{\Delta} })$. Thus
$$
a_{j+1} \le a_j \, \E \left\{ \exp\left( - c\, a_j \, Y\, \ee^{S_{\Delta-1}} \widetilde{M}_{j^3} \right) \right\} .
$$

\noindent As a consequence,
$$
a_{\lfloor n^{1/3}\rfloor} \le \prod_{j=0}^{\lfloor n^{1/3}\rfloor -1} \E \left\{ \exp\left( - c\, a_j \, Y\, \ee^{S_{\Delta-1}} \widetilde{M}_{j^3} \right) \right\} .
$$

We claim that for any collection of non-negative random variables $(\eta_j, \; 0\le j\le n)$ and $\lambda \ge 0$,
$$
\prod_{j=0}^n \E \left( \ee^{-\eta_j}\right) \le
\ee^{- \lambda} + \prod_{j=0}^n \P \{ \eta_j < \lambda\}.
$$

\noindent Indeed, without loss of generality, we can assume that $\eta_j$ are independent; then
$$
\prod_{j=0}^n \E \left( \ee^{-\eta_j}\right) \le \E \left( \ee^{- \max_{0\le j\le n} \eta_j} \right)
\le \ee^{- \lambda} + \P \left\{ \max_{0\le j\le n} \eta_j < \lambda \right\} = \ee^{- \lambda} + \prod_{j=0}^n \P \{ \eta_j < \lambda\} ,
$$

\noindent as claimed.

We have thus proved that
$$
a_{\lfloor n^{1/3}\rfloor} \le \ee^{-n} + \prod_{j=0}^{\lfloor n^{1/3}\rfloor-1} \P \left\{ c\, a_j \, Y\, \ee^{S_{\Delta-1}} \widetilde{M}_{j^3} < n\right\}.
$$

Recall that $a_j = \E (\varrho_{j^3+1})$. By the already proved lower bound in Theorem \ref{t:beta}, we have $a_j \ge \exp(-c_6\, j)$ for $j\ge j_0$. Hence, for $j_0\le j\le \lfloor n^{1/3}\rfloor -1$,
$$
\P \left\{ c\, a_j \, Y\, \ee^{S_{\Delta-1}} \widetilde{M}_{j^3} \ge n\right\} \ge \P \{ Y \ge 1\} \, \P \left\{ \widetilde{M}_{j^3} \ge {1\over j^3}\right\} \P \left\{ S_{\Delta-1} \ge c_6\, j + \log \left( {j^3 n\over c}\right) \right\} .
$$

\noindent Of course, $\P\{Y \ge 1\} = \ee^{-1}$; and by (\ref{M-ltail}), $\P \{ \widetilde{M}_{j^3} \ge {1\over j^3}\} = \P \{ M_{j^3} \ge {1\over j^3}\}\ge {1\over 2}$ for all large $j$. On the other hand,  since $\Delta -1\ge 3j^2$, we have $\P \{ S_{\Delta-1} \ge c_6\, j + \log ( {j^3 n\over c}) \} \ge c_{20}>0$ for large $n$ and all $j\ge \log n$. We have thus proved that, for large $n$ and some constant $c_{21}\in (0,\, 1)$,
$$
a_{\lfloor n^{1/3}\rfloor} \le \ee^{-n} + \prod_{j=\lceil \log n \rceil}^{\lfloor n^{1/3}\rfloor-1} (1-c_{21}) \le \exp \left( -c_{22} \, n^{1/3} \right) .
$$

\noindent Since $a_{\lfloor n^{1/3}\rfloor} = \E(\varrho_{\lfloor n^{1/3}\rfloor^3 +1}) \ge \E(\varrho_{n+1})$, this yields (\ref{E(beta)<<}), and thus the upper bound in Theorem \ref{t:beta}.\hfill$\Box$

\section{Proof of Theorem \ref{t:V*}: lower bound}
\label{s:V-lb}

Without loss of generality (see Remark \ref{r:V*}), we can assume $\psi'(1)=0$. In this case, the lower bound in Theorem \ref{t:V*} follows from the upper bound in Theorem \ref{t:beta} (proved in the previous section) by means of Proposition \ref{p:V*}, with $c_7:= c_3$.\hfill$\Box$

\section{Proof of Theorem \ref{t:main}}
\label{s:proof}

For the sake of clarity, Theorem \ref{t:main} is proved in two distinct parts.

\subsection{Upper bound}
\label{subs:ub}

We first assume $\psi'(1)=0$. By Theorem \ref{t:beta}, $P_\omega \{ \tau_n < \tau_0\} = \varrho_n \le \exp(- c_3 \, n^{1/3})$, $\P$-almost surely for all large $n$. Hence, by writing $L(\tau_n) := \# \{ 1\le i \le \tau_n: X_i =e \}$, we obtain: $\P$-almost surely for all large $n$ and any $j\ge 1$,
$$
P_\omega \{ L(\tau_n) \ge j \} = [ P_\omega \{ \tau_n > \tau_0 \} ]^j \ge [ 1- \ee^{- c_3 \, n^{1/3}} ]^j,
$$

\noindent which, by the Borel--Cantelli lemma, implies that, for any constant $c_{23}<c_3$ and $\p$-almost surely all sufficiently large $n$,
$$
L(\tau_n) \ge \ee^{c_{23} \, n^{1/3}}.
$$

\noindent Since $\{ L(\tau_n) \ge j \} \subset \{ X_{2j}^* < n\}$, we obtain the desired upper bound in Theorem \ref{t:main} (case $\psi'(1)=0$), with $c_2 := 1/(c_3)^3$.

To treat the case $\psi'(1)>0$, we first consider an RWRE $(Y_k, \; k\ge 0)$ on the half-line $\z_+$ with a reflecting barrier at the origin. We write $T_Y(y) := \inf\{ k\ge 0: \, Y_k =y\}$ for $y \in \z_+ \backslash \{ 0\}$. Then
$$
P_\omega \{ T_Y(y) \le m\} = \sum_{i=1}^m P_\omega \{ T_Y(y) =i \} \le \sum_{i=1}^m P_\omega \{ Y_i =y\} = \sum_{i=1}^m \omega^i (0,y),
$$

\noindent where, by an abuse of notation, we use $\omega(\cdot , \cdot)$ to denote also the transition matrix of $(Y_k)$. Since $(Y_k)$ is reversible, we have $\omega^i (0,y) = {\pi(y)\over \pi(0)} \omega^i (y,0)$, where $\pi$ is an invariant measure. Accordingly,
$$
P_\omega \{ T_Y(y) \le m\} \le \sum_{i=1}^m {\pi(y)\over \pi(0)} \omega^i (y,0) \le m {\pi(y)\over \pi (0)} .
$$

\noindent As a consequence, for any $n\ge 1$,
$$
P_\omega \{ T_Y(n) \le m\} \le \min_{1\le y \le n} P_\omega \{ T_Y(y) \le m\} \le m \min_{1\le y \le n} {\pi(y)\over \pi(0)} .
$$

\noindent It is easy to compute $\pi$: we can take $\pi(0)=1$ and
$$
\pi(y) := \sum_{z=1}^y \log {\omega(z,z-1)\over \omega(z,z+1)} , \qquad y\in \z_+ \backslash \{ 0\}.
$$

\noindent Therefore, for $n\ge 1$,
\begin{equation}
    P_\omega \{ T_Y (n) \le m\} \le m \min_{y\in\,
    ]\!] e, \, x]\!]} A(y) = m\, \ee^{- \overline{V}
    (x)},
    \label{diese}
\end{equation}

\noindent where $\overline{V}(x)$ is defined in (\ref{Vbar}).

We now come back to the study of $X$, the RWRE on $\T$. Fix $x\in \T_n$. Let $Z=(Z_k, \; k\ge 0)$ be the restriction of $X$ on the path $[\![ e, \, x]\!]$ as in (\ref{Xhat}). Let $T_Z(x) := \inf\{ k\ge 0: Z_k =x\}$. By (\ref{diese}), we have $P_\omega\{ T_Z(x) \le m\} \le m\, \ee^{- \overline{V}(x)}$. It follows from the trivial inequality $T(x) \ge T_Z(x)$ that
$$
P_\omega\{ \tau_n \le m\} \le \sum_{x\in \T_n} P_\omega\{ T(x) \le m\} \le \sum_{x\in \T_n} P_\omega\{ T_Z(x) \le m\} \le m \sum_{x\in \T_n}\ee^{- \overline{V}(x)} .
$$

Since $\psi'(1)>0$, we can consider $0<\theta<1$ as in (\ref{gamma}). Then
$$
\sum_{x\in \T_n} \ee^{- \overline{V}(x)} \le \exp\left( - (1-\theta)\min_{x\in \T_n} \overline{V}(x) \right) \sum_{x\in \T_n} \ee^{-\theta V(x)}.
$$

\noindent Since $\E(A^\theta)=1$, it is easily seen that $\sum_{x\in \T_n} \ee^{-\theta V(x)}$ is a positive martingale. In particular, $\sup_{n\ge 1} \sum_{x\in \T_n} \ee^{-\theta V(x)} <\infty$, $\P$-almost surely. On the other hand, according to Theorem \ref{t:V*}, we have $\min_{x\in \T_n} \overline{V}(x) \ge c_7 \, n^{1/3}$, $\P$-almost surely for all large $n$. Therefore, for any constant $c_{24}< (1- \theta) c_7$, we have
$$
\sum_n P_\omega \left\{ \tau_n \le \ee^{c_{24} \, n^{1/3}} \right\} < \infty, \qquadæ\hbox{\rm $\P$-a.s.,}
$$

\noindent from which the upper bound in Theorem \ref{t:main} (case $\psi'(1)>0$) follows readily, with $c_2 := 1/ [(1-\theta)c_7]^3$.\hfill$\Box$

\subsection{Lower bound}
\label{subs:lb}

By means of the Markov property, one can easily get a recurrence relation for $E_\omega (\tau_n)$, from which it follows that for $n\ge 1$,
\begin{equation}
    E_\omega(\tau_n) = {\gamma_n(e) \over
    \varrho_n},
    \label{E(tau)}
\end{equation}

\noindent where $\varrho_n$ and $\gamma_n(e)$ are defined by: $\beta_n(x) = 1$ and $\gamma_n(x) = 0$ (for $x\in \T_n$), and
\begin{eqnarray*}
    \beta_n(x)
 &=& {\sum_{i=1}^\deg A(x_i) \beta_n
    (x_i) \over 1+ \sum_{i=1}^\deg A(x_i) \beta_n
    (x_i)} ,
    \\
    \gamma_n(x)
 &=& {[1/\omega(x, \, {\buildrel \leftarrow \over
    x})] + \sum_{i=1}^\deg A(x_i) \gamma_n
    (x_i) \over 1+ \sum_{i=1}^\deg A(x_i) \beta_n
    (x_i)}, \qquad 1\le |x| \le n,
\end{eqnarray*}

\noindent and $\varrho_n := \sum_{i=1}^\deg \omega (e, e_i) \beta_n (e_i)$, $\gamma_n(e) := \sum_{i=1}^\deg \omega(e, e_i) \gamma_n (e_i)$. See Rozikov~\cite{rozikov} for more details. As a matter of fact, $\beta_n(x)$ (for $1\le |x| \le n$) is the same as the one introduced in (\ref{beta-recurrence}), and $\varrho_n$ can also be expressed as $P_\omega \{ \tau_n < \tau_0 \}$.

We claim that
\begin{equation}
    \sup_{n\ge 1} \, {\gamma_n(e)\over n} < \infty,
    \qquad \hbox{\rm $\P$-a.s.}
    \label{gamma<}
\end{equation}

By admitting (\ref{gamma<}) for the moment, we are able to prove the lower bound in Theorem \ref{t:main}. Indeed, in view of (the lower bound in) Theorem \ref{t:beta} and (\ref{E(tau)}), we have $E_\omega(\tau_n) \le c_{25}(\omega) \, n \exp( c_4 \, n^{1/3})$, $\P$-almost surely for all large $n$. It follows from Chebyshev's inequality and the Borel--Cantelli lemma that $\p$-almost surely for all sufficiently large $n$, $\tau_n \le c_{25}(\omega) \, n^3 \exp( c_4 \, n^{1/3})$, which yields
$$
\liminf_{n\to \infty} {X_n^*\over (\log n)^3} \ge {1\over (c_4)^3}, \qquad \hbox{$\p$-a.s.}
$$

\noindent This is the desired lower bound in Theorem \ref{t:main}.

It remains to prove (\ref{gamma<}). By the ellipticity condition, ${1\over \omega(x, \, {\buildrel \leftarrow \over x})} \le c_{26}$, so that
$$
\gamma_n(x) \le c_{26} + \sum_{i=1}^\deg A(x_i) \gamma_n (x_i) .
$$

\noindent Iterating the inequality, we obtain:
$$
\gamma_n(e) \le c_{26} \left( 1+ \sum_{j=1}^{n-1} \sum_{x\in \T_j} \prod_{y\in \, ] \! ] e_i, \, x] \! ]} A(y) \right) = c_{26} \left( 1+ \sum_{j=1}^{n-1} M_j \right), \qquad n\ge 2 .
$$

\noindent where $M_j$ is already introduced in (\ref{Mn}).

There exists $0<\theta\le 1$ such that $\E(A^\theta) ={1\over \deg}$: indeed, if $p={1\over \deg}$ and  $\psi'(1)=0$, then we simply take $\theta=1$, whereas if $p={1\over \deg}$ and $\psi'(1)>0$, then we take $0<\theta<1$ as in (\ref{gamma}). We have
$$
M_j^\theta \le \sum_{x\in \T_j} \prod_{y\in \, ] \! ] e_i, \, x] \! ]} A(y)^\theta.
$$

\noindent Since $j\mapsto \sum_{x\in \T_j} \prod_{y\in \, ] \! ] e_i, \, x] \! ]} A(y)^\theta$ is a positive martingale, we have $\sup_{j\ge 1} M_j <\infty$, $\P$-almost surely. This yields (\ref{gamma<}), and thus completes the proof of the lower bound in Theorem \ref{t:main}.\hfill$\Box$

\bigskip
\bigskip


{\footnotesize

\baselineskip=12pt

\noindent
\begin{tabular}{lll}
&\hskip20pt Yueyun Hu
    & \hskip45pt Zhan Shi\\
&\hskip20pt D\'epartement de Math\'ematiques
    & \hskip45pt Laboratoire de Probabilit\'es et
       Mod\`eles Al\'eatoires \\
&\hskip20pt Universit\'e Paris XIII
    & \hskip45pt Universit\'e Paris VI\\
&\hskip20pt 99 avenue J-B Cl\'ement
    & \hskip45pt 4 place Jussieu\\
&\hskip20pt F-93430 Villetaneuse
    & \hskip45pt F-75252 Paris Cedex 05\\
&\hskip20pt France
    & \hskip45pt France \\
&\hskip20pt {\tt yueyun@math.univ-paris13.fr}
    & \hskip45pt {\tt zhan@proba.jussieu.fr}
\end{tabular}

}

\end{document}